%% file: Switch_2.tex
\documentclass[5p]{elsarticle}

\usepackage{etoolbox}
\newtoggle{review}
\togglefalse{review}
\ifdim \baselinestretch pt = 1.5 pt 
  \toggletrue{review}
\fi

\usepackage{amssymb}

\usepackage{lineno}

\pdfoutput=1  
\usepackage{amsmath}
\allowdisplaybreaks[4]
\usepackage{algorithmic}
\usepackage{array}
\usepackage{color,soul}    
\usepackage[table]{xcolor} 
\usepackage{pifont}        
\usepackage{verbatim}
\usepackage{adjustbox}
\usepackage{comment}
\usepackage{booktabs}
\usepackage{makecell}
\usepackage{float}
\usepackage{multirow}
\usepackage{url}

\iftoggle{review}{
    \linenumbers
}{
    \usepackage[fontsize=10pt]{scrextend}
}

\journal{SoftwareX}

\begin{document}
\sloppy

\begin{frontmatter}



\title{Switch 2.0: A Modern Platform for Planning High-Renewable Power Systems}


\author[uhm,erg,jj]{Josiah Johnston\corref{cor}}
\ead{josiah.johnston@gmail.com}

\author[ucbee,pucc]{Rodrigo Henr\'iquez\corref{cor}}
\ead{rhenriqueza@uc.cl}

\author[pucc]{Benjam\'in Maluenda}
\ead{bmaluend@uc.cl}

\author[uhm]{Matthias Fripp\corref{pcor}}
\ead{mfripp@hawaii.edu}

\cortext[cor]{Corresponding author}
\cortext[pcor]{Principal corresponding author}

\address[uhm]{Department of Electrical Engineering, University of Hawaii at Manoa, Honolulu, Hawaii, USA}
\address[erg]{Energy and Resources Group, University of California, Berkeley, California, USA}
\address[jj]{Josiah Johnston, PhD Consulting \& Research}
\address[ucbee]{Department of Electrical Engineering and Computer Sciences, University of California, Berkeley, USA}
\address[pucc]{Department of Electrical Engineering, Pontificia Universidad Cat\'olica de Chile, Santiago, Chile}

\begin{abstract}
\input{Sections/abstract.tex}
\end{abstract}


\end{frontmatter}


\section{Motivation and significance}
\label{motivation}
\input{Sections/motivation.tex}

\section{Software description}
\label{software}
\input{Sections/software.tex}

\section{Illustrative Example}
\label{studies}
\input{Sections/case_studies.tex}



\section{Impact and Conclusions}
\label{conclusions}
\input{Sections/conclusion.tex}

\section*{Acknowledgements}
\label{acknowledgments}
The authors thank researchers who have used prior versions of Switch, whose discussions, collaborations and support have inspired continued work on this platform. Portions of the work reported here were funded by grants from the US Department of Transportation’s University Transportation Centers Program, Research and Innovative Technology Administration (PO\#291166), the National Science Foundation (\#1310634), and the Ulupono Initiative.



\bibliographystyle{elsarticle-num}
\bibliography{references}

\section*{Required Metadata}
\label{sec:metadata}

\section*{Current code version}
\label{sec:codeversion}

\begin{table}[!h]
\begin{tabular}{|l|p{6.5cm}|p{6.5cm}|}
\hline
\textbf{Nr.} & \textbf{Code metadata description} & \textbf{Please fill in this column} \\
\hline
C1 & Current code version & 2.0.0 \\
\hline
C2 & Permanent link to code/repository used for this code version & https://github.com/switch-model/switch/releases/tag/2.0.0 \\
\hline
C3 & Legal Code License   & Apache-2.0 \\
\hline
C4 & Code versioning system used & git \\
\hline
C5 & Software code languages, tools, and services used & python, Pyomo \\
\hline
C6 & Compilation requirements, operating environments & python 2.7 \\
\hline
C7 & If available Link to developer documentation/manual & http://switch-model.org \\
\hline
C8 & Support email for questions & Matthias Fripp \textless{}mfripp@hawaii.edu\textgreater{} \\
\hline
\end{tabular}
\caption{Code metadata (mandatory)}
\label{tab:metadata}
\end{table}

\section*{Current executable software version}
\label{tab:codeversion}

\begin{table}[!h]
\begin{tabular}{|l|p{6.5cm}|p{6.5cm}|}
\hline
\textbf{Nr.} & \textbf{(Executable) software metadata description} & \textbf{Please fill in this column} \\
\hline
S1 & Current software version & 2.0.0 \\
\hline
S2 & Permanent link to executables of this version  & https://pypi.org/project/switch-model/2.0.0/ \\
\hline
S3 & Legal Software License & Apache-2.0 \\
\hline
S4 & Computing platforms/Operating Systems & BSD, Linux, OS X, macOS, Microsoft Windows, Unix-like \\
\hline
S5 & Installation requirements & python 2.7; pip or git \\
\hline
S6 & If available, link to user manual - if formally published include a reference to the publication in the reference list & http://switch-model.org \\
\hline
S7 & Support email for questions & Matthias Fripp \textless{}mfripp@hawaii.edu\textgreater{}\\
\hline
\end{tabular}
\caption{Software metadata (optional)}
\label{tab:optionalmetadata}
\end{table}

\end{document}

%% file: Sections/abstract.tex
Switch 2.0 is an open-source platform for planning transitions to low-emission electric power systems, designed to satisfy 21st century grid planning requirements. Switch is capable of long- and short-term planning of investments and operations with conventional or smart grids, integrating large shares of renewable power, storage and/or demand response. Applications include integrated resource planning, investment planning, economic and policy analyses as well as basic research. Switch 2.0 includes a complete suite of power system elements, including unit commitment, part-load efficiency, planning and operating reserves, fuel supply curves, storage, hydroelectric networks, policy constraints and demand response. It uses a robust, modular architecture that allows users to compose models by choosing built-in modules \textit{\`a la carte} or writing custom modules.

%% file: Sections/motivation.tex

Climate change mitigation requires electric power systems to achieve deep emission cuts by mid-century while serving increasingly electrified transport and heating sectors \cite{williams2012technology, wei2013deep}. Power system infrastructure lasts for many decades, so capacity planners must already consider this transition. Meanwhile, the cost of renewable power and storage are falling below fossil fuels \cite{XcelEnergy2017,IRENA2018,Kittner2017} and ubiquitous computation is enabling ``smart grids'' to reschedule demand to help balance intermittent wind and solar power \cite{OConnell2014, Paterakis2017}. These changes create a planning challenge of unprecedented complexity: how to select the optimal portfolio of intermittent and controllable power sources, storage, transmission and demand response, to move toward 100\% renewable power over the coming decades, while maintaining an exact balance between supply and demand every hour?

To meet this challenge, power system stakeholders---utilities, transmission operators, regulators, environmental and consumer organizations and academic researchers---need general-purpose power system planning software with five interrelated capabilities: multiple investment steps over several decades, sequential modeling of individual hours of operation, detailed modeling of unit commitment decisions for generators, user-extendability, and low-cost or open-source licensing. In this paper, we introduce Switch 2.0, the first power system model to include all of these elements in a single package, enabling a broad range of users to confidently plan 21st century power systems.

\input{Figures/bigtable.tex}

Below, we discuss each of these capabilities and the contributions of previous software, summarized in Table \ref{bigtable}. We cannot include a complete list of useful features or models, so we focus on an important subset. Note that each of the capabilities is provided by \emph{some} previous software, but Switch 2.0 is the first to provide \emph{all} of them simultaneously. This new fusion is essential, because omitting any of these elements significantly impairs planning for high-renewable power systems.

\emph{Multiple investment periods} and \emph{Inter-hour relationships}. Co-optimizing capacity investment decisions for several periods over a multi-decade horizon---rather than making a single-step investment plan---reduces the risk of creating stranded assets in an evolving grid \cite{Pfenninger:2014ht}. Further, operational decisions should be modeled during chronological sequences of timesteps within many study days. This allows direct modeling of choices that depend on the sequence of hours within a single day, such as unit commitment (slow and expensive power plant start/stop decisions), when to charge and discharge finite storage, or when to schedule time-shiftable demand within the day \cite{Welsch:2015kk, nweke2012benefits, poncelet2016impact, Pfenninger:2014ht, Palmintier15, wogrin-system-states, de-jonghe-dr, denholm2015overgeneration, rosenkranz2016analyzing, palmintier2011impact}. Modeling these accurately will become more important as the power system moves toward 100\% renewables and resources adapt intensively to each day's weather conditions. Meeting both of these requirements is computationally difficult. A number of capacity planning models use multiple investment periods but don't model sequential  timepoints within each day, impairing their ability to study unit commitment, storage and demand response \cite{loulou2004documentation,eia2016NEMS,short2011regional,heaps2016long,Karlsson:2008gc,Shawhan:2014eg,gil2015hydro,nweke2012benefits}; others use a single planning step with hour-by-hour timeseries  \cite{oemof_developer_group_2017,dorfner_urbs_2018,BrownPyPSAPythonPower2018,Palmintier15,stiphout2017impact,o2013model,jenkins_enhanced_2017}. The PLEXOS model in long-term planning mode \cite{nweke2012benefits} uses hourly resolution for peak days, but only 1-2 blocks on all other days, obscuring the year-round view. Version 1 of Switch pioneered the combination of multiple investment periods and inter-hour behavior in a single model, by using multiple, separate sample timeseries per planning period \cite{fripp2012switch}. Resolve \cite{RESOLVE2017}, based on Switch 1, uses the same approach, and two newer models---RPM \cite{mai_resource_2013} and WIS:dom \cite{macdonald_future_2016}---use a similar approach. Switch 2.0 introduces a more general, flexible version of this approach.

\emph{Generator unit commitment}. Power plant unit commitment---decisions about the ``lumpy'' and slow start/stop behavior of conventional power plants, and the inefficiency of operating at part-load---significantly restrict power system flexibility and drive the cost of providing operating reserves \cite{Padhy2004a}. This is increasingly important as power systems ramp conventional plants more intensively to counterbalance intermittent renewable power \cite{nweke2012benefits, poncelet2016impact, Pfenninger:2014ht, Welsch:2015kk, jonghe2011deter, stiphout2017impact, palmintier2011impact, Palmintier15}. Unit commitment is omitted or simplified in many long-term models \cite{loulou2004documentation,eia2016NEMS,short2011regional,heaps2016long,Karlsson:2008gc,Shawhan:2014eg,fripp2012switch,RESOLVE2017}. However it must be modeled directly for some pressing questions, e.g., to find the optimal balance between batteries, demand response and partly loaded generators to provide reserves for uncertain renewable power. Several new planning models include direct modeling of unit commitment and operating reserves \cite{young_us-regen_2018, oemof_developer_group_2017, dorfner_urbs_2018, BrownPyPSAPythonPower2018, Palmintier15, stiphout2017impact, o2013model, jenkins_enhanced_2017, fripp2012switch, RESOLVE2017, mai_resource_2013, macdonald_future_2016}. (PLEXOS LT Plan also includes these, but is limited by coarse time sampling \cite{nweke2012benefits}.) Switch 2.0 adds full unit commitment capabilities to the Switch framework.

\emph{User-defined extensions}. For many studies, users need to add new technologies or subsystems to planning models, e.g., to study hydrogen networks, advanced demand response, or heat and power distribution systems. Some existing models have a formal framework for user extensions \cite{loulou2004documentation, oemof_developer_group_2017, BrownPyPSAPythonPower2018}. Others are not available to the public, making customizability moot \cite{short2011regional, young_us-regen_2018, jenkins_enhanced_2017, mai_resource_2013}, or allow some extensions but not fundamental changes to the model \cite{exemplar2016plexos}. Open-source models allow free revision \cite{heaps2016long, Karlsson:2008gc, Shawhan:2014eg, dorfner_urbs_2018, fripp2012switch, RESOLVE2017}, but often only by creating new forks of the whole model for each change. Syncing these forks and sending changes back to the main model is difficult, impairing innovation and reuse. One of the strongest features of Switch 2.0 is an exceptionally flexible, modular architecture that allows users to arbitrarily extend the model, share new features between projects, and easily contribute new elements back to the main codebase. This framework also allows plug-and-play selection among built-in modules, so users can reconfigure models with no new code at all.

\emph{Low-cost or open software}. Low-cost or open-source planning models allow small research groups, nonprofits or residents of less wealthy countries to evaluate new technologies, policies and plans, ensuring that more voices are heard and better solutions are found. Open-source software supports replication of research findings and agreement among stakeholders \cite{Wilson:2014ck, Pfenninger:2014ht, DeCarolis:2012dg, atatech:2016}. A number of planning models are only available for use by their authors and partners (although model formulations are sometimes available to the public) \cite{short2011regional,young_us-regen_2018,jenkins_enhanced_2017,mai_resource_2013,Palmintier15,stiphout2017impact,o2013model}. Commercial models may cost more than \$25,000/year \cite{exemplar2016plexos}, excluding less-wealthy stakeholders. Several models are released as open-source software that runs within lower-cost commercial optimization environments, making them more accessible \cite{loulou2004documentation, eia2016NEMS, Karlsson:2008gc, Shawhan:2014eg}. Switch 2.0 joins a small group of open-source models with no proprietary dependencies \cite{heaps2016long, oemof_developer_group_2017, dorfner_urbs_2018, BrownPyPSAPythonPower2018, RESOLVE2017}, allowing use by the broadest possible range of users.

%% file: Figures/bigtable.tex

\newcommand*\rcell[1]{\rotatebox{90}{\makecell[l]{#1}}}
\newcommand*\OK{\cellcolor{black!25}\ding{51}}

\newcommand*\notreleased{-}
\newcommand*\no{}
\newcommand*\partly{-}
\newcommand*\notproprietary{$\dagger$}
\newcommand*\framework{*}
\newcommand*\nosource{?}
\newcommand*\plexosblocks{a}
\newcommand*\gasonly{b}
\newcommand*\bidsonly{c}
\newcommand*\daylong{d}
\newcommand*\experimental{e}

\begin{table}[ht!]
\centering
\caption{Feature comparison among leading electricity capacity planning models}
\iftoggle{review}{
    \begin{adjustbox}{width=0.75\textwidth}
}{
    \begin{adjustbox}{width=0.485\textwidth}
}
\label{bigtable}
\setlength{\tabcolsep}{0.2em} 
\begin{tabular}{|c|c|c|c|c|c|}\hline  
    {Model} & \rcell{Multiple \\investment periods$\,$} & \rcell{Inter-hour \\relationships} & \rcell{Generator \\unit commitment$\,$} & \rcell{User-defined \\extensions} & \rcell{Low-cost or \\open software} \\\hline
    \multicolumn{1}{|l|}{TIMES/MARKAL \cite{loulou2004documentation,contaldi2007evaluation}} & \OK & \no & - & \OK & - \\\hline
    \multicolumn{1}{|l|}{NEMS \cite{eia2016NEMS,eia2012national,gabriel2001national,eia_availability_2018}} & \OK & \no & - & \no & - \\\hline
    \multicolumn{1}{|l|}{ReEDS \cite{short2011regional}} & \OK & \no & - & n/a & \no \\\hline
    \multicolumn{1}{|l|}{LEAP + OseMOSYS \cite{heaps2016long,howells2011osemosys}} & \OK & \no & - & - & \OK \\\hline
    \multicolumn{1}{|l|}{Balmorel \cite{Karlsson:2008gc}} & \OK & \no & - & - & - \\\hline
    \multicolumn{1}{|l|}{E4 Simulation Tool \cite{Shawhan:2014eg}} & \OK & \no & \no & - & - \\\hline
    \multicolumn{1}{|l|}{US-REGEN (non-UC mode) \cite{young_us-regen_2018}} & \OK & \no & \OK & n/a & \no \\\hline
    \multicolumn{1}{|l|}{US-REGEN (UC mode) \cite{young_us-regen_2018}} & \no & \OK & \OK & n/a & \no \\\hline
    \multicolumn{1}{|l|}{Oemof \cite{oemof_developer_group_2017,Hilpert:2017eq}} & \no & \OK & \no & \OK & \OK \\\hline
    \multicolumn{1}{|l|}{URBS \cite{dorfner_urbs_2018,huber2015optimizing}} & \no & \OK & \no & - & \OK \\\hline
    \multicolumn{1}{|l|}{PyPSA \cite{BrownPyPSAPythonPower2018}} & \no & \OK & \OK & \OK & \OK \\\hline
    \multicolumn{1}{|l|}{Palmintier and Webster \cite{Palmintier15}} & \no & \OK & \OK & - & \no \\\hline
    \multicolumn{1}{|l|}{Stiphout et al. \cite{stiphout2017impact}} & \no & \OK & \OK & - & \no \\\hline
    \multicolumn{1}{|l|}{O'Neill et al. \cite{o2013model}} & \no & \OK & \OK & - & \no \\\hline
    \multicolumn{1}{|l|}{GenX \cite{jenkins_enhanced_2017}} & \no & \OK & \OK & n/a & \no \\\hline
    \multicolumn{1}{|l|}{PLEXOS LT Plan \cite{exemplar2016plexos,gil2015hydro,nweke2012benefits}} & \OK & - & - & - & \no \\\hline
    \multicolumn{1}{|l|}{Switch 1.x \cite{fripp2012switch, nelson2012high}} & \OK & \OK & - & - & - \\\hline
    \multicolumn{1}{|l|}{RESOLVE \cite{RESOLVE2017, CAPUC2017IRP}} & \OK & \OK & - & - & \OK \\\hline
    \multicolumn{1}{|l|}{RPM \cite{mai_resource_2013,nrel_resource_2018}} & \OK & \OK & \OK & n/a & \no \\\hline
    \multicolumn{1}{|l|}{NEWS/WIS:dom \cite{macdonald_future_2016,clack_weather-informed_2018}} & \OK & \OK & ? & ? & \no \\\hline
    \multicolumn{1}{|l|}{\textbf{Switch 2.0 }} & \OK & \OK & \OK & \OK & \OK \\\hline

\multicolumn{6}{l}{\makecell[l]{
\ding{51} fully supported \\
- partially supported (see text for details) \\
n/a model is not available to outside researchers \\
? information not available
}}
\end{tabular}
\end{adjustbox}
\end{table}


%% file: Sections/software.tex
\subsection{Software Architecture}
\label{architecture}
Switch 2.0 is a Python package that can be installed via standard Python package managers or directly from its Github repository. Installation instructions are at \url{http://switch-model.org}. Switch uses the open-source Pyomo \cite{hart2017pyomo} optimization framework to define models, load data and solve instances. Models can be solved using any optimization software compatible with Pyomo, which includes most commercial and open-source solvers.

\begin{figure*}[t]
\centering
\includegraphics[width=.99\textwidth]{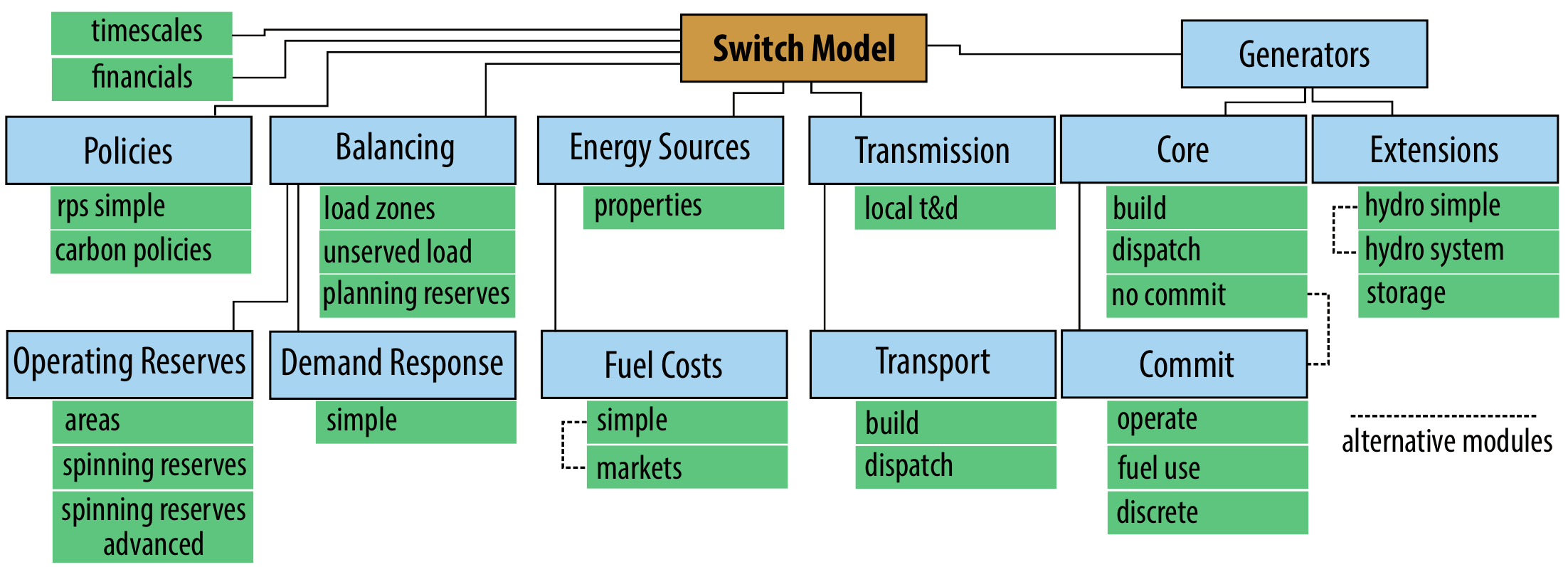}
\caption{Package and module structure of Switch 2.0. Blue boxes are subpackages, green boxes are modules.}
\label{fig:modules}
\end{figure*}

Switch 2.0 uses a fine-grained, modular approach to define power system models, which allows the formulation to be easily customized for the needs of each study. This modular architecture reflects the modularity of actual power systems, where individual elements operate independently but contribute to the system's total costs and power balance. Core modules in Switch define spatially and temporally resolved  balancing constraints for energy and reserves, and an overall system cost. Separate modules represent components such as generators, batteries or transmission links. These modules interact with the overall optimization model by adding terms to the shared energy and reserve balances and the overall cost expression. They can also define additional decision variables and constraints to govern operation of each technology or subsystem. This allows technologies to be packaged in plug-and-play modules that participate as fully integrated components of the overall model.

Each built-in or user-supplied Switch module is implemented by creating a Python module file that defines one or more standard callback functions that will be called at each stage of generating and solving a model: defining and parsing command-line arguments, defining model components, defining costs or energy balance equations that are shared between modules, loading data from an input directory,
and performing post-solve functions.
These functions are detailed in the Supplementary Material.

Users configure the model by creating a text file containing a list of modules to be used. At runtime, Switch loads each module and runs through definition, compilation, solution, and export stages, calling callback functions of each module in turn, starting with core modules to define the basic framework, followed by specialized modules to define custom technologies or policies.

This system is highly flexible, making it easy to add or subtract from the codebase, typically without having to modify the built-in modules. Basic extensions are simple to write, and advanced modules have an unlimited ability to extend the model.
By changing the choice of modules, users can also switch easily between distinct modeling modes, such as running a sparse capacity-expansion model, followed by a more detailed production-cost assessment of the proposed portfolio.

\subsection{Software Functionalities}
\label{functionalities}

Figure \ref{fig:modules} presents the subpackages and modules that are included in Switch 2.0. The key modules and subpackages are summarized below. A complete mathematical formulation of the model is provided in the Supplementary Material, and thoroughly documented source code may be inspected in the model repository.

\emph{Timescales}.
This module defines three timescales for decision making: \emph{periods} of one or more years where investment decisions are made, \emph{timepoints} within each period when operational decisions are made, and \emph{timeseries} that group timepoints into chronological sequences. Timepoints within each timeseries have a fixed duration specified in hours, and timeseries have a fixed weight that denotes how many times this type of series is expected to occur in the corresponding period. This approach can represent any standard time structure: a load duration curve (many one-hour timeseries per planning period), a collection of sample days during each period (several one-day timeseries), or an 8760-hour timeseries as typically used by production cost models (a single, year-long timeseries).

\emph{Financials}.
This module defines the objective function and financial parameters for the model. Other modules may register investment or operational cost components with this module. Switch then minimizes the net present value of all costs over the entire study.

\emph{Balancing}.
This subpackage defines \emph{load zones}, geographic regions with load timeseries in which energy supply and demand must be balanced in all timepoints. Other modules may register power injections or withdrawals with the power balance constraints. An optional \emph{Unserved Load} module allows imbalances with a user-specified penalty per unit of energy. This package also includes subpackages and modules for \emph{Planning Reserves}, \emph{Operating Reserves} and \emph{Demand Response}, which are described in more detail in the Supplemental Material.

\emph{Generators}.
This subpackage defines all possible generation projects. The \emph{Core} subpackage describes construction and operation constraints and decisions for all projects standard thermal generators or renewables without storage. A simple economic dispatch may be chosen via the \emph{No Commit} module, or a full unit commitment formulation can be incorporated through the \emph{Commit} subpackage. Generators may be assigned a single energy source or allowed to switch optimally between fuels in order to meet targets for emissions or renewables.

Additional power sources may be implemented via optional modules in the \emph{Extensions} subpackage within the \emph{Generators} package. The \emph{Storage} module defines a generic framework for storage technologies, such as pumped hydro, batteries, flywheels, and others. This  formulation permits independent sizing of energy and power elements. The \emph{Hydro System} module represents a cascading water network operating in parallel with the energy network, whereas the \textit{Hydro Simple} module merely enforces average water availability for each sampled timeseries. Switch 2.0 offers both linearized unit commitment and mixed-integer unit-commitment via the \emph{Commit} subpackage and \emph{Commit.Discrete} module. When using the \emph{Commit} subpackage, users may also optionally provide multi-segment heat rate curves, minimum up and down times and startup costs and energy.

Switch aggregates generators into generation ``projects''. These are stacks of one or more similar generating units in the same transmission zone, but not necessarily at the same site. Capacity can be added to each generation project in different years, and then portions of the available capacity are committed or dispatched as needed. This approach significantly reduces model size. 

\emph{Transmission}.
Switch offers several approaches for representing transmission network capabilities. The most basic strategy is to use a single-zone (\textit{copperplate}) formulation, which ignores restrictions on spatial transfer of power. Alternatively, a \textit{transport model} can be used to represent transmission capabilities in a simplified manner \cite{fripp2012switch,loulou2004documentation, schaber_transmission_2012}. In a transport model, the study region is divided into zones which are internally well-connected but have constrained connections (flowgates) to neighboring zones. The size of each flowgate and the amount of power transfer are decision variables. Transport models are designed to approximate the capabilities of the network and the cost of improvements without modeling the electrical behavior of the network directly. They provide an attractive balance between granularity and tractability in expansion models, because they use only linear terms, even when network expansion is considered. Switch uses a copperplate model by default, or a transport model can be adopted by using the \emph{Transport} subpackage of the \emph{Transmission} subpackage.


The optional \emph{Local T\&D} module within the \emph{Transmission} subpackage represents power transfers from the zonal node to customers, as well as distributed energy resources. This module enables a simplified consideration of the impact of distributed generation, efficiency or demand response on distribution network investments or losses.

\emph{Energy Sources}
This subpackage defines fuel and non-fuel energy sources. Fuel costs can be either represented by a \emph{Simple} flat cost per period or through a \emph{Markets} module which supports supply curves and regional markets which may or may not be interconnected. These are important for power systems that may need to switch to new energy sources such as biofuels or LNG in order to meet strict targets for renewables or emissions

\emph{Policies}
This subpackage defines investment and operational policies. Current modules include enforcing a \emph{Simple RPS} (Renewable Portfolio Standard) and \emph{Carbon Policies}, such as carbon taxes and caps. Investment or production cost credits can be modeled by adjusting fixed or variable costs.



\emph{Testing}.
Switch uses an automated testing framework to support quality assurance and control as users develop code. The testing framework includes a mix of unit, integration, and regression tests.
This helps ensure that changes to the model formulation don't introduce bugs that break existing models, allowing faster and more robust evolution of the codebase.

%% file: Sections/case_studies.tex

In 2015, the State of Hawaii adopted a 100\% renewable power target by 2045. Here we use Switch 2.0 and datasets prepared for a forthcoming study of that transition, to assess the potential benefits of obtaining operating reserves from several different types of battery and/or demand-side response. This is an important question for engineers and policymakers planning the transition to a 100\% renewable power system, and can be addressed directly using the combination of long-term modeling, chronological daily modeling, unit commitment, operating reserves, storage and customizability offered in Switch 2.0. We are not aware of any other published models that can address these questions directly.

This study uses assumptions from the Hawaiian Electric Company's latest integrated resource plan \cite{HECOHawaiianElectricCompanies2016b}; model configuration and data are available from ref. \cite{HawaiiReservesStudyRepository}. There are four types of lithium-ion batteries: contingency-oriented batteries can make 10 deep cycles per year and regulation-oriented batteries can make up to 15,000 shallow cycles per year. Neither of these can provide bulk load-shifting from hour to hour. However, load-shifting batteries can provide 4 or 6 hours of energy storage and complete up to 365 cycles per year. Peak demand could be reduced up to 10\% via demand response.

In this case study, we use Switch 2.0 to answer several questions about these resources.
(1) If load-shifting batteries could provide contingency and/or regulating reserves while charging or discharging, how much money would that save?
(2) How much money could be saved per household by implementing demand response as a simple load-shifting service?
(3) How much more could be saved if demand response also provided contingency and regulating reserves?


To address these questions, we ran Switch using these standard modules:
\begin{itemize}
\itemsep-0.4em   
\item \verb|timescales|,
\item \verb|financials|,
\item \verb|balancing.load_zones|,
\item \verb|energy_sources.properties|,
\item \verb|generators.core.build|,
\item \verb|generators.core.dispatch|,
\item \verb|reporting|,
\item \verb|energy_sources.fuel_costs.markets|,
\item \verb|generators.core.proj_discrete_build|,
\item \verb|generators.core.commit.operate|,
\item \verb|generators.core.commit.fuel_use|,
\item \verb|generators.core.commit.discrete|,
\item \verb|generators.extensions.storage|,
\item \verb|balancing.operating_reserves.areas|, and
\item \verb|balancing.operating_reserves|\linebreak[0]{}\verb|.spinning_reserves_advanced|.
\end{itemize}

These comprise Switch's core formulation plus the following elements: discrete-sized generators, detailed unit-commitment, two-way flow for batteries, regional fuel markets, and spinning reserve targets. The formulation of these modules and the Hawaii reserve rule are detailed in the Supplementary Material. We also used several modules from a \verb|hawaii| subpackage, to define options for fuel market expansion, demand-response, electric-vehicle charging and operating rules for some individual generators.

We used Switch's scenario-solving system to define five scenarios: (1) "battery bulk": load-shifting batteries only provide bulk inter-hour load-shifting, not reserves; there is no demand response (DR), and electric vehicles (EVs) charge at business-as-usual times; (2) "battery bulk and conting": same as "battery bulk", but load-shifting batteries can also provide contingency reserves; (3) "battery bulk and reg": same as "battery bulk and conting", but load-shifting batteries can also provide regulating reserves; (4) "DR bulk": same as "battery bulk and reg" plus DR can provide bulk load-shifting (i.e., the model can move up to 10\% of demand from each hour to any other hour, provided it doesn't raise demand by more than 80\% in any hour), and EVs charge at optimal times each day; no reserves from DR or EVs; (5) "DR bulk and reg": same as "DR bulk", but DR and EVs can also provide up and down contingency and regulation reserves equal to the difference between the amount of load scheduled each hour and the minimum and maximum allowed loads.

We then used Switch to select the optimal investment plan for each of the five scenarios, with investment decisions made every five years from 2020 through 2045. For these scenarios, we considered 12 sample days during each period, using weather from the 23rd day of each month in 2008. Costs are scaled as if each sample day were repeated 152 times, filling out the five-year period.

\begin{figure*}[t]
\centering
\includegraphics[width=7in]{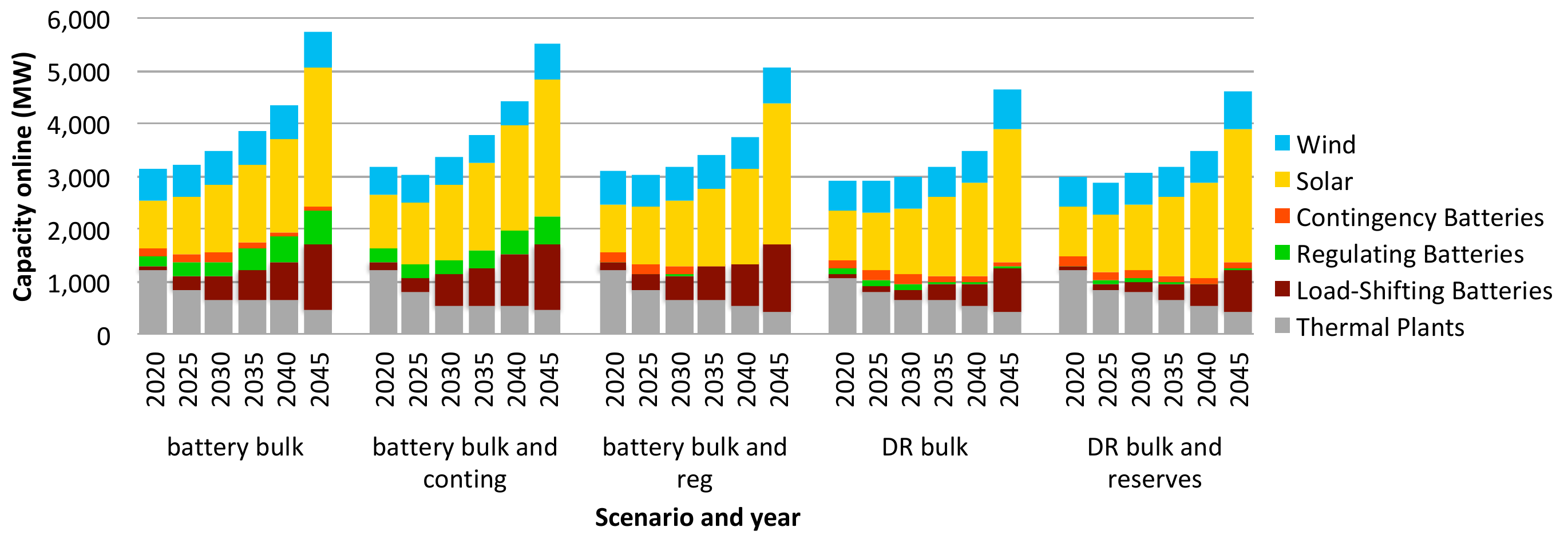}
\caption{Generation and storage capacities in five scenarios.}
\label{fig:hawaii_capacity_by_scenario}
\end{figure*}

Figure \ref{fig:hawaii_capacity_by_scenario} shows the optimal generation and storage portfolios selected in the five scenarios.
From this brief study, we can draw several useful conclusions:

(1) Obtaining contingency reserves from load-shifting batteries is not likely to provide large savings (only about \$225 per customer on an NPV basis); this is because contingency batteries make up only a small share of the system's assets, and are relatively inexpensive per MW of capacity.

(2) Obtaining regulating reserves from load-shifting batteries is financially attractive (additional \$615 savings per customer). These savings occur because this control strategy reduces the need to build dedicated batteries for regulation.

(3) Inter-hour load shifting via demand response (e.g., in response to real-time pricing) could save an additional \$1,850 per customer, an attractive option.

(4) Providing regulation and contingency reserves from demand response saves only an additional \$159 per customer. The savings come mainly from a small reduction in the need for regulating batteries. This may not be cost-effective, unless this service can ``piggyback'' on the price-based response infrastructure.

%% file: Sections/conclusion.tex


Switch 2.0 is an open-source platform that can perform power system studies such as investment planning, production cost simulations, or economic and policy analyses. Its modular architecture allows users to formulate and solve models with customized features and varying levels of complexity. Software best practices such as emphasizing readability of code, automated testing, embedded documentation, and review processes allow effective collaborative extension of modules. All of these characteristics make Switch an ideal platform for researchers, educators, students, nonprofit organizations and members of industry to study evolving power systems.

Switch 1.0 and 2.0 have been used by researchers at several universities and nonprofits to analyze the evolution of power systems in many regions and countries. These include the Western Electrical Coordinating Council (western portion of the United States and Canada with part of Mexico), Chile, Hawaii, Mexico, China, Nicaragua, Japan, Kenya, East Africa, and Peru \cite{fripp2012switch, nelson2012high, mileva2013sunshot, wakeyama_impact_2015, wei2013deep, PoncedeLeonBarido:2015ec, Sanchez:2015ij, he2016switch}.

More recently, several projects have advanced in tandem with Switch 2.0, and improvements from those projects have been incorporated into the Switch 2.0 codebase. Here we give three examples of active areas of research with Switch 2.0. Each of these studies is made possible by the unique combination of temporal scale, operational detail, customizability and transparency that Switch 2.0 offers.

\emph{Renewable Expansion in Chile}.
Switch 2.0 is being used to investigate optimal expansion of renewables in Chile's relatively high-hydro power system. The study spans from 2020 to 2040 and includes 20 load zones, 23 transmission corridors, and 107 generation projects. Some are existing projects that may be expanded with new units and some are proposed projects that have no existing capacity. Hydroelectric projects are either located at a water reservoir or act as run-of-river plants (see the Supplementary Material for an illustration). This hydraulic network modeling structure was implemented in \cite{benjamin2017paper} as a user-defined module, and has since been integrated into the main code repository.

\emph{Hawaii Renewable Energy Planning}.
In 2015, Hawaii adopted legislation requiring its electric utilities to reach 100\% renewable power by 2045. Switch has taken a central role in planning to meet this target. Hawaiian Electric Company (HECO) used RESOLVE, based on Switch 1.x, to develop its first integrated resource plan after the RPS was adopted \cite{HECOHawaiianElectricCompanies2016b}. In parallel with this, stakeholders used Switch 2.0 to evaluate and critique this plan and offer alternatives \cite{fripp_making_2016, heco_hawaiian_2016, FrippIncentiveProblemsPlanning2017}. This work led to consensus that it would be cost-effective to add several hundred megawatts of additional wind and solar to the company's 5-year plan (in a 1200 MW power system), significantly accelerating the State's move from high-cost oil to lower-cost renewables.

\emph{Demand Response}.
Switch 2.0 is especially well-suited to studying the use of price-responsive demand to help integrate renewable power. To support this, Switch 2.0 includes an advanced demand response module that can find equilibrium between the power system design and any convex demand system. This has been used for studies of the economic benefits of real-time pricing in a high-renewable power system \cite{ImeldaVariablePricing2018, ImeldaVariablePricingSocialforthcoming}, and the effect of improved timing of electric vehicle charging on the design and operation of high-renewable power systems \cite{DasSavingsPeakReduction2015}.

\emph{Validation}. Switch 2.0's production-cost capability has been benchmarked against the industry-standard GE MAPS model and found to agree closely for 18 renewable adoption scenarios in Oahu and Maui \cite{fripp_inter-model_nodate}.
\newline
\newline
\emph{Future work}. Future work will focus on parallel solutions across large regions and many timesteps, modules for DC and AC power flow, more definitions of operational reserves, case studies of stochastic programming, and preparation of public, national datasets. Efforts are also underway to expand the research community built around Switch and foster international collaboration and discussion. Switch-Mexico researchers have also been prototyping a graphical web interface \cite{switchGUI2016}.